\documentclass[12pt]{article}
\newcommand{\lap}{\mbox{$\bigtriangleup$}}

\newcommand{\ra}{{\mbox{$\rightarrow$}}}
\newcommand{\be}{\begin{equation}}
\newcommand{\ee}{\end{equation}}
\newtheorem{mrem}{Remark}
\newtheorem{mthm}{Theorem}
\newtheorem{mpro}{Proposition}

\newtheorem{thm}{Theorem}[section]

\newtheorem{lem}{Lemma}[section]

\newtheorem{rem}{Remark}[section]
\begin{document}

\title{Super Polyharmonic Property of Solutions for PDE Systems and Its Applications}
\author{Wenxiong Chen
\thanks{Partially supported by NSF Grant DMS-0604638}
\hspace{.2in}
Congming Li
\thanks{Partially supported by NSF Grants DMS-0908097 and EAR-0934647}}
\date{}
\maketitle

\begin{abstract}
In this paper, we prove that all the positive solutions for the PDE system
\begin{equation}
(- \lap )^k u_i = f_i(u_1, \cdots, u_m) , \;\; x \in R^n , \; i = 1, 2, \cdots, m
\label{PDEsA}
\end{equation}
are super polyharmonic, i.e.
$$ (- \lap)^j u_i > 0 , \;\; j=1, 2, \cdots, k-1; \; i =1, 2, \cdots, m.$$

To prove this important super polyharmonic property, we introduced a few new ideas and
derived some new estimates.

As an interesting application, we establish the equivalence between the integral system
\begin{equation}
u_i(x) = \int_{R^n} \frac{1}{|x-y|^{n-\alpha}} f_i(u_1(y), \cdots, u_m(y)) d y , \;\; x \in R^n
\label{intsaa}
\end{equation}
and PDE system (\ref{PDEsA}) when $\alpha = 2k<n.$

In the last few years, a series of results on qualitative properties for solutions of
integral systems (\ref{intsaa}) have been obtained, since the introduction of a powerful tool--
the method of moving planes in integral forms. Now due to the equivalence established here, all these properties can be applied to the
corresponding PDE systems.

We say that systems (\ref{PDEsA}) and (\ref{intsaa}) are equivalent, if whenever
$u$ is a positive solution of (\ref{intsaa}), then $u$ is also a solution of
$$(- \lap )^k u_i =  c f_i(u_1, \cdots, u_m) , \;\; x \in R^n , \; i= 1,2, \cdots, m $$
with some constant $c$; and vice versa.
\end{abstract}

\medskip


\section{Introduction}

The {\em method of moving planes in integral forms} was introduced in \cite{CLO} to solve an open problem posed by Lieb \cite{Lie}: Can one classify all the positive solutions of the integral equation
\begin{equation}
u(x) = \int_{R^n} \frac{1}{|x-y|^{n-\alpha}} u^{\frac{n+\alpha}{n-\alpha}}(y) d y , \;\; x \in R^n \; ?
\label{1.1}
\end{equation}
This is an Euler equation of the functional associated with the well-known Hardy-Littlewood-Sobolev inequality in a special case. The classification would provide the best constant in the corresponding inequality.

To show the radial symmetry of solutions for PDEs, usually a method of moving planes is applied, which relies heavily on maximum principles. While integral equations do not possess such local properties. To carry on the moving of planes, we explore global features of the integral equations and estimate certain integral norms. This is the essence of the {\em method of moving planes in integral forms}. It can be very conveniently applied to (\ref{1.1}) to establish radial symmetry of the solutions under quite mild integrability assumptions, that is $u \in L^{\frac{n+\alpha}{n-\alpha}}_{loc}$. As compare to the corresponding PDE
\begin{equation}
(- \lap)^{\alpha/2} u = u^{\frac{n+\alpha}{n-\alpha}} , \;\; x \in R^n ,
\label{de}
\end{equation}
to carry on the method of moving planes, one at least need to require the solution be $\alpha$ times differentiable and can only do this when $\alpha$ is an even number. While fractional powers of Laplacian in (\ref{de}) have many applications in mathematical physics and
related fields, such as anomalous diffusion and quasi-geostrophic
flows, turbulence models and water waves, molecular dynamics and relativistic quantum mechanics
of stars \cite{BoGe} \cite{CaVa} \cite{Co} \cite{MaMcTa} \cite{TaZa}, as well as in probability
and finance \cite{A} \cite{Be} \cite{CT}.

 Since its introduction, the {\em method of moving planes in integral forms} has become a powerful tool to establish symmetry, a priori estimates, and non-existence for a variety of integral equations and systems, including
\smallskip

(i) a system of $m$ equations \cite{CL2}
\begin{equation}
\left\{\begin{array}{l} u_i(x)  =
\displaystyle{\int_{R^n} \frac{1}{|x-y|^{n-\alpha}}}
f_i(u(y))dy, \;\; x \in {R^n} , \;\; i=1,...m,
 \\
 0< \alpha <n, \; \mbox{ and } \; u(x)=(u_1(x), u_2(x), \cdots, u_m(x));
\end{array}
\right.
\label{ELLS}
\end{equation}

(ii) the fully nonlinear integral systems involving Wolff potentials \cite{CL5}:
\begin{equation}
\left\{\begin{array}{ll}
u(x) = W_{\beta, \gamma}(v^q)(x) , & x \in R^n; \\
v(x) = W_{\beta, \gamma} (u^p)(x) , & x \in R^n;
\end{array}
\right.
\label{WFP}
\end{equation}
where
$$ W_{\beta,\gamma} (f)(x) = \int_0^{\infty} \left[ \frac{\int_{B_t(x)} f(y) dy}{t^{n-\beta\gamma}} \right]^{\frac{1}{\gamma-1}} \frac{d t}{t};$$

(iii) equations and systems involving Bessel potentials \cite{MC1} \cite{MC2}:
$$
\left\{\begin{array}{ll}
u(x) = g_{\alpha} \ast f(u,v) , \; & x \in R^n , \\
v(x) = g_{\beta} \ast g(u,v) , \; & x \in R^n ,
\end{array}
\right.
$$
where $\ast$ is the convolution and
$$ g_{\alpha} (x) = \frac{1}{(4\pi)^{\alpha} \Gamma (\frac{\alpha}{2})} \int_0^{\infty} \exp \left\{
- \frac{\pi |x|^2}{t} - \frac{t}{4\pi} \right\} \frac{d t}{t^{(n-\alpha)/2 +1}} $$
is the kernel of the Bessel potential,

 (iv) boundary values problems on half spaces \cite{LZ} and in unit balls \cite{CZ},

 (v) problems on complete Riemannian manifolds \cite{QR}.
 \smallskip

There are many more applications of this method to integral systems, and the interested readers may see \cite{CL1} \cite{CL3} \cite{CL6} \cite{CLO1} \cite{CLO2}  \cite{Ha} \cite{HWY} \cite{JL} \cite{Liu} \cite{LiY} \cite{LM} \cite{LQ} \cite{LSL} \cite{MC} \cite{MCL} \cite{MZ} and the references therein.

In \cite{JL2} \cite{LL} \cite{LLM}, asymptotic analysis for solutions of integral systems near the origin and at infinity were also obtained.

If one can show that an integral system is equivalent to a PDE system, then the results obtained for the integral system can be carried over to the PDE system.

For instance, in \cite{CLO}, we showed that (\ref{1.1}) is equivalent to semi-linear elliptic PDE (\ref{de}) and thus classified all the positive solutions for both integral equation (\ref{1.1}) and PDE (\ref{de}).

In \cite{QR},
Qing and Raske  considered a smooth family of equations
\begin{equation}
P_{\alpha}[g] u = u^{\frac{n+\alpha}{n-\alpha}}
\label{QR1}
\end{equation}
on Riemannian manifolds.
When $\alpha = 2$, $P_{\alpha}[g]$ is the conformal Laplacian, and when $\alpha =4$, it is the Paneitz operator. They showed that, under certain conditions, (\ref{QR1}) is equivalent to a conformally invariant integral equation. Using the {\em method of moving planes in integral forms}, they obtained a priori estimate for the positive solutions for the integral equation as well as for PDE (\ref{QR1}); and based on this estimate and a degree theory, they proved the existence of solutions.

In \cite{MZ}, Ma and Zhao were able to classified the positive solutions for the stationary Choquard equation
\begin{equation}
\lap u - u + 2 u \cdot \left( \frac{1}{|x|} \ast |u|^2 \right) = 0 , \;\;  u \in H^1(R^3)
\label{Cho}
\end{equation}
and hence answered an open question posed by Lieb \cite{Lie1}. Their idea was to use Bessel potential to write equation (\ref{Cho}) as an equivalent integral system
\begin{equation}
\left\{\begin{array}{ll}
u = g_2 \ast (u v) \\
v = \frac{2}{|x|} \ast |u|^2 .
\end{array}
\right.
\label{Cho1}
\end{equation}
Then they applied the {\em method of moving planes in integral forms} to derive the radial symmetry of the positive solutions of (\ref{Cho1}), and hence of (\ref{Cho}) due to the equivalence.

Since the {\em method of moving planes in integral forms} is applied directly to integral equations and systems, whether or not these results can be extended to the corresponding PDEs depends on whether one can establish the equivalence between the integral equations and PDEs. Is the integral system (\ref{ELLS}) equivalent to a PDE system? What can we say about  (\ref{WFP})?

To establish the equivalence, a general approach is to multiply both sides of the PDEs by the Green's function on the ball of radius $R$, integrating by parts, then letting $R \ra \infty$ and taking limits. In order to show that the boundary terms tend to zero, one important ingredient is to obtain the super polyharmonic property of the solutions for the PDEs. This is the main objective of our present paper.

Consider the following more general system of PDE inequalities:

\begin{equation}
(- \lap )^{k} u_i \geq f_i(u) , \;\; x \in R^n , \; i = 1, 2, \cdots, m.
\label{sysa}
\end{equation}
where $u = (u_1, \cdots, u_m)$. We prove

\begin{mthm} Let $u=(u_1, \cdots, u_m)$ be a positive solution of (\ref{sysa}).
Assume that
\begin{equation}  f_i(u) \geq 0 , \;\; i = 1, 2, \cdots, m .
\label{f1}
\end{equation}
Let $w = u_1 + u_2 + \cdots + u_m$. Suppose there exist $p>1$, $\delta >0$, and $C, C_{\delta} >0$, such that
\begin{equation}
\sum_{i=1}^m f_i(u) \geq \left\{\begin{array}{ll} C_{\delta} w^p & \mbox{ for $w$ sufficiently large}\\
C_{\delta} & \mbox{ if } w \geq \delta .
\end{array}
\right.
\label{f2}
\end{equation}

Then we have
\begin{equation}
(- \lap)^j u_i > 0 , \;\; x \in R^n, \;\; j = 1, 2, \cdots k-1; \;\; i=1, 2, \cdots, m .
\label{4.1}
\end{equation}
\label{mthmsysm}

\end{mthm}

When $f_i(u) = u^{p_i}$, the power of $- \lap$ need not to be the same, as illustrated in the following
system of two inequalities:

\begin{equation}
\left\{\begin{array}{ll}
(- \lap )^t u \geq v^q (x) , & x \in R^n \\
(- \lap)^s v \geq u^p (x) , & x \in R^n
\end{array}
\right.
\label{1.3}
\end{equation}
with positive integers $t$ and $s$.

\begin{mthm} Assume that $p, q>1$. Let $(u,v)$ be a pair of positive solutions for inequality system (\ref{1.3}).
Then
$$ \left\{\begin{array}{ll}
(-\lap)^i u > 0 , & x \in R^n, \;\; i=1, 2, \cdots, t-1;\\
(-\lap)^i v > 0 , & x \in R^n, \;\; i=1,2, \cdots, s-1.
 \end{array}
 \right.
 $$
\label{mthmsys2}
\end{mthm}

\begin{mrem} One can verify that when (and only when) $s=t$, Theorem \ref{mthmsys2} is a special case of Theorem \ref{mthmsysm}.
\end{mrem}

Theorem \ref{mthmsysm} contains the following result of Lin \cite{Lin} (on a fourth order equation) and Wei and Xu \cite{WX} (on general even order equations) as a special case:

\begin{mpro} Let $u$ be a positive solution of
$$ (-\lap)^k u = u^p (x), \; \; x \in R^n $$
with $p>1$, then
$$ (-\lap)^i u > 0 , \; i = 1, 2, \cdots, k-1. $$
\label{mprowx}
\end{mpro}

To prove Theorem \ref{mthmsysm}, we introduce some new ideas and obtained a few interesting estimates.

It is well-known that the super polyharmonic property (\ref{4.1}) is very useful in analyzing the corresponding system. Now, as an important and immediate application, we use Theorem \ref{mthmsysm} and \ref{mthmsys2} to establish the equivalence between integral systems and the corresponding PDE systems.

We say that the integral system
\begin{equation}
u_i(x) = \int_{R^n} \frac{1}{|x-y|^{n-\alpha}} f_i(u(y)) d y , \;\; x \in R^n, \; i=1,2, \cdots, m
\label{intsa}
\end{equation}
and the PDE system
\begin{equation}
(- \lap )^{\alpha/2} u_i =  f_i(u) , \;\; x \in R^n , \; i =1,2, \cdots, m
\label{PDEsc}
\end{equation}
 are equivalent, if whenever
$u=(u_1, \cdots, u_m)$ is a positive solution of (\ref{intsa}), then $u$ is also a solution of
$$(- \lap )^{\alpha/2} u_i =  c f_i(u) , \;\; x \in R^n , \; i =1,2, \cdots, m$$
with some constant $c$; and vice versa.
\medskip

We prove

\begin{mthm} Let $t$ and $s$ be positive integers less than $\frac{n}{2}$, and assume that $p, q >1$. Then
the PDE system
\begin{equation}
\left\{\begin{array}{ll}
(- \lap )^t u = v^q (x) , & x \in R^n \\
(- \lap)^s v = u^p (x) , & x \in R^n
\end{array}
\right.
\label{1.3a}
\end{equation}
is equivalent to the integral system
$$\left\{\begin{array}{ll}
u(x) = \int_{R^n} \frac{1}{|x-y|^{n-2t}} v^q(y) d y , & x \in R^n, \\
v(x) = \int_{R^n} \frac{1}{|x-y|^{n-2s}} u^p(y) dy , & x \in R^n.
\end{array}
\right.
$$
\label{mthmequiv2}
\end{mthm}

\begin{mthm}
Let $\alpha$ be an even integer less than $n$. Then under the assumption of Theorem \ref{mthmsysm},  PDE system (\ref{PDEsc}) is equivalent to integral system (\ref{intsa}).
 \label{mthmequiv}
\end{mthm}

For other real values of $\alpha$, we define the
positive solutions of
\begin{equation}
(- \lap )^{\alpha/2} u_i =   f_i(u) , \;\; x \in R^n , \;\; i=1, \cdots, m
\label{PDEsa}
\end{equation}
 in the distribution sense, i.e. $u
\in H^{\frac{\alpha}{2}}(R^n)$ and satisfies
\be
\int_{R^n} (- \lap)^{\frac{\alpha}{4}} u_i \, (- \lap)^{\frac{\alpha}{4}} \phi \, dx =
\int_{R^n} f_i(u(x)) \, \phi (x) \, dx , \;\; i=1, \cdots, m;
\label{6.20a}
\ee
for any
$\phi \in C_0^{\infty} (R^n)$ and $\phi (x) \geq 0.$ Here, as usual,
$$\int_{R^n} (- \lap)^{\frac{\alpha}{4}} u_i \, (-
\lap)^{\frac{\alpha}{4}} \phi \, dx$$
is defined by Fourier transform
$$\int_{R^n} |\xi|^{\alpha} \hat{u_i} (\xi) \overline{\hat{\phi} (\xi)} \,
d\xi , $$ where $\hat{u_i}$ and $\hat{\phi}$ are the Fourier
transform of $u_i$ and $\phi$ respectively.

\begin{mthm}
PDE system (\ref{PDEsa}) as defined above is
equivalent to integral system (\ref{intsa}).
\label{mthmequivalpha}
\end{mthm}

Based on these
equivalences, a series of previous results obtained by many authors for integral equations and systems can now be applied to the corresponding PDE systems.
\medskip

In order to better illustrate our new idea, in Section 2, we will apply it to a single inequality.  In Section 3, we apply our new idea to a system of two inequalities and prove Theorem \ref{mthmsys2}. In Section 4, we consider a more general system of $m$ inequalities and obtain Theorem \ref{mthmsysm}. Finally, in Section 5, we establish the equivalences between integral and PDE systems and thus prove Theorem \ref{mthmequiv2}, \ref{mthmequiv}, and \ref{mthmequivalpha}.
\medskip

We would like to mention that, in Section 4 and 5, when considering a general system, we introduce an interesting idea, so that we can reduce the treatment of a system into that of a single equation or inequality and hence simplify the proofs remarkably.

Throughout this paper, as usual, when we say ``a solution'' to a PDE or inequality system, we mean ``a classical solution'' if not otherwise indicate.

\section{Super Polyharmonic Properties for Single Inequalities}

 Let $p$ be a positive integer and $q>1$. Consider
 \begin{equation}
 (- \lap)^p u \geq u^q(x) , \;\; x \in R^n .
 \label{1}
 \end{equation}

 \begin{thm} For each positive solution $u$ of (\ref{1}), it holds
 \begin{equation}
(- \lap u )^i > 0 , \;\; i = 1, \cdots, p-1.
\label{2}
\end{equation}
\label{thmWX}
\end{thm}

\leftline{\bf Proof of Theorem \ref{thmWX}.}

Write
 $$ v_i = (- \lap)^i u , \;\; i = 1, \cdots, p-1. $$

{\bf {\em Part I.} }

We first show that
\begin{equation}
v_{p-1} (x) > 0 .
\label{3}
\end{equation}

Suppose in contrary, there are two possible cases:

Case i) There exists $x^o \in R^n$, such that
$$ v_{p-1} (x^o) < 0 .$$

Case ii) $v_{p-1}(x) \geq 0$ and there is a point $\tilde{x}$, such that
$$ v_{p-1}(\tilde{x}) = 0 .$$
In this case, $\tilde{x}$ is a local minimum of $v_{p-1}$, and we must have $- \lap v_{p-1} (\tilde{x}) \leq 0$. This contradict with $$ - \lap v_{p-1} = u^q > 0 .$$

Therefore we only need to consider Case i).
Without loss of generality, we may assume that $x^o = 0.$

{\em Step 1.}

In this step, we will show that, after a few times of {\em re-centers}, if we denote the
resulting functions by $\tilde{v}_k$, then
\begin{equation}
\tilde{v}_{p-2}> 0, \; \tilde{v}_{p-3} < 0, \; \tilde{v}_{p-4} >0, \cdots .
\label{B1}
\end{equation}
The positiveness of $u$ implies that $p$ must even. Hence by this alternating sign nature, we have
\begin{equation}
- \lap \tilde{u} \equiv \tilde{v}_1 < 0 ,
\label{B2}
\end{equation}
and therefore
\begin{equation}
\tilde{u}(r) \geq \tilde{u}(0) > 0 .
\label{B3}
\end{equation}
Here, we assume that after {\em re-centers}, $\tilde{u}$ is radially symmetric about the origin.
\medskip

To derive (\ref{B1}),
let
$$ \bar{u} (r) = \frac{1}{|\partial B_r (0)|} \int_{\partial B_r (0)} u \, d \sigma  $$
be the average of $u$. Then by Jensen's inequality and the well-known property that
$$ \overline{\lap u} = \lap \bar{u}, $$we have
\begin{equation}
\left\{\begin{array}{ll}
- \lap \bar{v}_{p-1} \geq \bar{u}^q ,\\
- \lap \bar{v}_{p-2} = \bar{v}_{p-1} ,\\
\cdots \cdots \\
- \lap \bar{u} = \bar{v}_1 .
\end{array}
\right.
\label{z1}
\end{equation}

From the first inequality, we have
\begin{equation}
\bar{v}_{p-1}'(r) < 0 \;\; \mbox{ and } \;\; \bar{v}_{p-1} (r) \leq \bar{v}_{p-1} (0) < 0 , \;\; \forall \; r >0 .
\label{4}
\end{equation}

Then from the second equation, we have
\begin{equation}
 - \frac{1}{r^{n-1}} \left( r^{n-1} \bar{v}'_{p-2} \right)' = \bar{v}_{p-1} (r) < \bar{v}_{p-1} (0) \equiv - c .
 \label{5}
 \end{equation}
Integrating yields
$$ \bar{v}_{p-2}'(r) > c r  \;  \; \mbox{ and } \; \bar{v}_{p-2}(r) \geq \bar{v}_{p-2}(0) + c r^2 , \;\; \forall \; r \geq 0 .$$
Hence there exists an $r_o$, such that $\bar{v}_{p-2} (r_o) > 0$.  Take a point $x^o$ with $|x^o| = r_o$ as the new center to
make average, then
$$ \bar{\bar{v}}_{p-2} (0) > 0 .$$
Obviously, $\bar{\bar{u}}$ and $\bar{\bar{v}}_i, \, i=1, \cdots, p-1,$ still satisfy (\ref{z1}).
Applying (\ref{4}) to $\bar{\bar{v}}_{p-1}$ instead of $\bar{v}_{p-1}$, we still have
\begin{equation}
\bar{\bar{v}}_{p-1} (r) < 0 , \; \forall \, r \geq 0 .
\label{6}
\end{equation}
From (\ref{5}) and repeat the same argument as for $\bar{v}_{p-2}$, we obtain
\begin{equation}
\bar{\bar{v}}_{p-2} (r) > 0 , \; \forall \, r \geq 0.
\label{7}
\end{equation}

Continuing this way, after a few steps of {\em re-centers} (denote the resulting functions by $\tilde{v}_k$ and $\tilde{u}$), we arrive at
\begin{equation}
(-1)^i \tilde{v}_{p-i} (r) > 0 , \; \forall \, r \geq 0 .
\label{8}
\end{equation}

This implies immediately that $p$ must be even, since if $p$ is odd, then (\ref{8}) implies that $u$ should be negative
somewhere, a contradiction with our assumption that $u>0$. Hence in the following, we assume that $p$ is even.

Let $$ u_{\lambda} (x) = \lambda^{2p/(q-1)} u(\lambda x) $$
be the re-scaling of $u$. Then equation (\ref{1}) is invariant under this re-scaling, and for any $\lambda >0$, $u_{\lambda}$ is still a positive solution of (\ref{1}). By (\ref{8}),
$$ - \lap \tilde{u} < 0 ,$$
and hence
$$ \tilde{u}'(r) > 0 .$$
It follows that
$$ \tilde{u}(r) > \tilde{u}(0) = c_o > 0 .$$
Hence one can choose $\lambda$ sufficiently large, such that $u_{\lambda}$ be as large as we wish. Without loss of generality,
we may assume that for any given $a_o > 0$, we already have
\begin{equation}
\tilde{u}(r) \geq a_o \geq a_o r^{\sigma_o} , \;\; \forall \, 0 \leq r \leq 1 .
\label{9}
\end{equation}
Here we choose $\sigma_o$, such that
\begin{equation}
\sigma_o \geq 1 \;\; \mbox{ and } \; \sigma_o q \geq 2p+n .
\label{10}
\end{equation}
The value of $a_o$ will be determined later.

From (\ref{9}),
$$ - \frac{1}{r^{n-1}} \left( r^{n-1} \tilde{v}'_{p-1} \right)' \geq a_o^q r^{\sigma_o q} .$$
Integrating both sides twice and taking into account of (\ref{8}) yield
$$ \tilde{v}_{p-1}(r) \leq - \frac{a_o^q}{(\sigma_o q +n)(\sigma_o q +2)}r^{\sigma_o q +2} .$$
This implies
$$ - \frac{1}{r^{n-1}} \left( r^{n-1} \tilde{v}'_{p-2} \right)' \leq - \frac{a_o^q}{(\sigma_o q +n)(\sigma_o q +2)}r^{\sigma_o q +2}.$$
And consequently
$$ \tilde{v}_{p-2}(r) \geq  \frac{a_o^q}{(\sigma_o q +n)(\sigma_o q +2)(\sigma_o q +n+2)(\sigma_o q +4)}r^{\sigma_o q +4}.$$

Continuing this way, we arrive at
\begin{equation}
\tilde{u}(r) \geq \frac{a_o^q}{(\sigma_o q +n+2p)^{2p}} r^{\sigma_o q +2p} \geq \frac{a_o^q}{(\sigma_o q + m )^m} r^{\sigma_o q + m} .
\label{11}
\end{equation}
Here we have denoted $2p+n$ as $m$.

Denote
$$ \sigma_{k+1} = 2 \sigma_k q \geq \sigma_k q + m \; \mbox{ and } \; a_{k+1} = \frac{a_k^q}{(2\sigma_k q)^m}.$$
Then obviously,
\begin{equation}
\tilde{u}(r) \geq a_1 r^{\sigma_1} .
\label{12}
\end{equation}

Suppose we have
$$ \tilde{u}(r) \geq a_k r^{\sigma_k} .$$
Then go through the entire process as above, we obtain
\begin{equation}
\tilde{u}(r) \geq a_{k+1} r^{\sigma_{k+1}} .
\label{13}
\end{equation}

Choose $l$ and then $\sigma_o$, such that
\begin{equation}
l(q-1)>2 \;\; \mbox{ and } \;\; \sigma_o \geq 2^{l+q+1} \, q^{2(l+1)+q} .
\label{14}
\end{equation}
We want to use induction to show that
\begin{equation}
a_k^q \geq (\sigma_k q)^{m(l+1)} , \;\; k = 0, 1, 2, \cdots .
\label{15}
\end{equation}

Obviously, we can make (\ref{15}) true for $k=0$ by choosing $a_0$ sufficiently large. Assume it is true
for $k$. Then we have
\begin{eqnarray*}
\frac{a_{k+1}^q}{(\sigma_{k+1} q)^{m(l+1)}} & =& \frac{\left[\frac{a_k^q}{(2\sigma_k q)^m} \right]^q}{(\sigma_{k+1} q)^{m(l+1)}} \\
&\geq& \left[\frac{\sigma_k}{2^{l+q+1} q^{2(l+1)+q}} \right]^m \geq 1.
\end{eqnarray*}
Therefore, (\ref{15}) is true for all integer $k$.

Now by (\ref{12}) and (\ref{15}),
$$ \tilde{u}(1) \geq a_k \geq (\sigma_k q)^{m(l+1)/q} \ra \infty , \;\; \mbox{ as } k \ra \infty .$$
This is obviously impossible. Therefore (\ref{3}) must hold.
\medskip

{\bf {\em Part II.}}

Base on the positiveness of $v_{p-1}$, we can now show that all other $v_k$ must also be positive:
$$ v_{p-i} (x) > 0 , \; \mbox{ for } i = 2, 3, \cdots, p-1.$$

Suppose for some $i$, $v_{p-i}$ is negative somewhere, then through the same arguments as in Step 1 of Part I,
after a few steps of {\em re-centers}, if we denote the resulting functions by $\tilde{v}_k$, then
the signs of $\tilde{v}_k$ are alternating, and by the positiveness of $\tilde{u}$, we must have
\begin{equation}
- \lap \tilde{u} < 0 \; \mbox{ and } \; \tilde{u} > 0 .
\label{B4}
\end{equation}
It follows that
$$\tilde{u} \geq c_o > 0 .$$

Coming back to the original equation
$$ - \lap v_{p-1} = u^q $$
we derive
$$ - \lap \tilde{v}_{p-1} \geq \tilde{u}^q \geq c_1 > 0 .$$
By a direct integration as in Part I, we would arrive at
$$ \tilde{v}_{p-1} < 0 , \; \mbox{ somewhere }.$$
This contradicts with the proven fact that $v_{p-1} > 0$, and thus completes the proof.

\section{A System of Two Inequalities}

Consider a system of two inequalities
in $R^n$
\begin{equation}
\left\{\begin{array}{ll}
(- \lap )^t u \geq v^q (x) ,\\
(- \lap)^s v \geq u^p (x) ,
\end{array}
\right.
\label{16}
\end{equation}
with positive integers $t$ and $s$. An entirely similar approach will be applied to systems of more inequalities and more general right hand sides functions in the next section. Nevertheless, the system in the next section does not include (\ref{16}), because there we require the powers of all $-\lap$ be the same, and here in (\ref{16}), they are allowed to be different.

Let
$$u_{\lambda}(x) = \lambda^{\alpha} u(\lambda x) \;\; \mbox{ and } \;\; v_{\lambda} (x) = \lambda^{\beta} v(\lambda x) , $$
where
$$ \alpha = \frac{2(t+sq)}{pq-1} \;\; \mbox{ and } \;\; \beta = \frac{2(s+tp)}{pq-1} .$$
Then one can easily verified that system (\ref{16}) is invariant under this re-scaling.

Let $$u_k = (-\lap)^k u \;\; \mbox{ and } \;\; v_k = (-\lap)^k v .$$
We want to show that
\begin{equation}
u_k , v_k > 0 \;\; \mbox{ for } k = 1, 2, \cdots \mbox{ up to $t-1$ or $s-1$ }.
\label{17}
\end{equation}

{\em {\bf Part I.}}

We first show that
\begin{equation}
u_{t-1} \equiv (-\lap)^{t-1} u > 0 .
\label{18}
\end{equation}

Suppose in contrary, there is some point $x^o$, such that
$$ u_{t-1} (x^o) < 0.$$
Then use the same argument as we did in deriving (\ref{8}) for the single equation, we obtain
\begin{equation}
(-1)^i \tilde{u}_{t-i} (r) > 0 \;\; \forall \, r \geq 0, \;\; i = 1, 2, \cdots, t-1 ,
\label{19}
\end{equation}
where $\tilde{u}_k$ is obtained from $u_k$ through several {\em re-centers}. (\ref{19}) implies that $t$ must be even.

From (\ref{19}),
$$ -\lap \tilde{u} (r) < 0 , $$
and hence
$$ \tilde{u}'(r) > 0 \;\; \mbox{ and } \; \tilde{u}(r) > \tilde{u}(0) > 0 .$$
From this and the inequality
$$ - \lap \tilde{v}_{s-1} \geq \tilde{u}^p (r) $$
we derive that $\tilde{v}_{s-1}$ must be negative somewhere. Now repeat the above argument and after a few steps of
{\em re-centers}, we arrive at
\begin{equation}
(-1)^i \tilde{v}_{s-i} (r) > 0 , \;\; \forall \; r \geq 0, \;\; i = 1, 2, \cdots, s-1 .
\label{20}
\end{equation}
Here we still denote the resulting functions after {\em re-centers} by $\tilde{v}_k$.
From here, we also obtain
$$ \tilde{v} (r) > \tilde{v}(0) > 0 .$$
After re-scaling, we can make $\tilde{v}_{\lambda}$ as large as we wish. Hence, without loss of generality, we may
assume that, for any fixed $a_0>0$,
\begin{equation}
\tilde{v} (r) \geq a_o \geq a_o r^{\sigma_o} , \;\; \forall \; 0 \leq r \leq 1 .
\label{21}
\end{equation}

Using this and $$ (-\lap)^t \tilde{u} \geq \tilde{v}^q $$
and through the same argument as in the previous section, we arrive at
\begin{equation}
\tilde{u}(r) \geq \frac{a_o^q}{(\sigma_o q +n+2t)^{2t}} r^{\sigma_o q +2t} \geq \frac{a_o^q}{(\sigma_o q + m )^m} r^{\sigma_o q + m} .
\label{22}
\end{equation}
Here we have let $m = \max\{n+2t, n+2s\}$.

Choose $\sigma_o$ large so that $\sigma_o q \geq m$.
Then (\ref{22}) becomes
\begin{equation}
\tilde{u} (r) \geq \frac{a_o^q}{(2\sigma_o q)^m} r^{2\sigma_o q} \equiv {b_o} r^{\eta_o}.
\label{23}
\end{equation}
Using the second equation
and repeat the same process as in deriving (\ref{22}), we have
\begin{equation}
\tilde{v}(r) \geq \frac{b_o^p}{(2\eta_o p)^m} r^{2\eta_o p} = \frac{a_o^{pq}}{2^{m(p+2)} q^{m(p+1)} p^m \sigma_o^{m(p+1)}} r^{4\sigma_o pq}
\equiv \frac{a_o^h}{c \sigma_o^{m(p+1)}} r^{4\sigma_o h}.
\label{24}
\end{equation}

let
$$ a_{k+1} = \frac{a_k^h}{c \sigma_k^{m(p+1)}} \;\; \mbox{ and } \;\; \sigma_{k+1} = 4\sigma_k h , \;\; k = 0, 1, 2, \cdots.$$
Then similar as in the single equation case, we can show by induction that
\begin{equation}
 a_k \ra \infty \;\; \mbox{ as } \; k \ra \infty .
 \label{25}
 \end{equation}
 Instead of carrying out the detail, we rather do it heuristically, which will better illustrate the idea.

 First notice that $\sigma_k \ra \infty$. To derive (\ref{25}), it suffice to show that
 \begin{equation}
 a_k^h \geq c \sigma_k^{m(p+1)+l} ,
 \label{26}
 \end{equation}
 for some $l$ to be determined later. We again use induction. We can choose $a_o$ large, so (\ref{26}) is true for $k=0$.
 Now assume (\ref{26}) holds for $k$, and we want to show that it is true for $k+1$. In fact, we have
 \begin{equation}
 \frac{a_{k+1}^h}{c \sigma_{k+1}^{m(p+1)+l}} \geq \frac{\sigma_k^{l(h-1)-m(p+1)}}{c (4h)^{m(p+1)+l}} .
 \label{27}
 \end{equation}

 Now choose $l$, such that the power $l(h-1)-m(p+1)$ of $\sigma_k$ is positive, say greater that $1$, and also choose $\sigma_o$
 large (note $\sigma_k \geq \sigma_o$) to ensure that (\ref{26}) holds for $k+1$. This verifies (\ref{18}). Similarly, one can derive that
 $$ v_{s-1} (r) \equiv (- \lap)^{s-1} v > 0 .$$

 {\bf {\em Part II.}}

Base on the positiveness of $v_{s-1}$, we can now show that all other $u_k$ must also be positive:
\begin{equation}
 u_{t-i} (x) > 0 , \; \mbox{ for } i = 2, 3, \cdots, t-1.
 \label{C1}
 \end{equation}

Suppose for some $i$, $u_{t-i}$ is negative somewhere, then through the same arguments as in Part I,
after a few steps of re-centers, if we denote the resulting functions by $\tilde{u}_k$, then
the signs of $\tilde{u}_k$ are alternating, and by the positiveness of $\tilde{u}$, we must have
\begin{equation}
- \lap \tilde{u} < 0 \; \mbox{ and } \; \tilde{u} > 0 .
\label{B4}
\end{equation}
It follows that
$$\tilde{u} \geq c_o > 0 .$$

Coming back to the original inequality
$$ - \lap v_{s-1} \geq u^q $$
we derive
$$ - \lap \tilde{v}_{s-1} \geq \tilde{u}^q \geq c_1 > 0 .$$
By a direct integration as in Part I, we would arrive at
$$ \tilde{v}_{s-1} < 0 , \; \mbox{ somewhere }.$$
This contradicts with the proven fact that $v_{s-1} > 0$ and thus (\ref{C1}) must hold.
Similarly, based on the positiveness of $u_{t-1}$ we can derive
$$ v_{s-i} (x) > 0 , \; \mbox{ for } i = 2, 3, \cdots, s-1.$$

This completes the proof.

\section{More General Systems}

Now we consider a system of more general inequalities
\begin{equation}
(- \lap )^{k} u_i \geq f_i(u) , \;\; x \in R^n , \; i = 1, 2, \cdots, m.
\label{sys}
\end{equation}
where $u = (u_1, \cdots, u_m)$.  We prove
\begin{thm} Let $u=(u_1, \cdots, u_m)$ be a positive solution of (\ref{sys}).
Assume that
\begin{equation}  f_i(u) \geq 0 , \;\; i = 1, 2, \cdots, m .
\label{f1}
\end{equation}
Let $w = u_1 + u_2 + \cdots + u_m$. Suppose there exist $p>1$, $\delta >0$, and $C, C_{\delta} >0$, such that
\begin{equation}
\sum_{i=1}^m f_i(u) \geq \left\{\begin{array}{ll} C_{\delta} w^p & \mbox{ for $w$ sufficiently large}\\
C_{\delta} & \mbox{ if } w \geq \delta .
\end{array}
\right.
\label{f2}
\end{equation}

Then we have
\begin{equation}
(- \lap)^j u_i > 0 , \;\; j = 1, 2, \cdots k; \;\; i=1, 2, \cdots, m .
\label{4.1}
\end{equation}
\label{thm4.1}

\end{thm}

{\bf Proof.}

Here we introduce an interesting idea, so that we can reduce the problem for system into one for single equation.

Suppose (\ref{4.1}) is not true. Without loss of generality, we may assume for some $j_o < k$,
\begin{equation}
(- \lap )^{j_o} u_1 < 0 , \;\; \mbox{ somewhere. }
\label{4.2}
\end{equation}

Let $$w_{\epsilon} = u_1 + \epsilon (u_2 + \cdots + u_m) $$ for some $\epsilon >0.$ Then by (\ref{f2}), we have
\begin{equation}
(- \lap )^k w_{\epsilon} \geq f_1 (u) + \epsilon (f_2(u) + \cdots + f_m(u)) \geq \epsilon C_{\delta} w_{\epsilon}^p .
\label{4.3}
\end{equation}

Now for each $\epsilon >0$, applying the same arguments as for single equation in Section 2, we can derive
$$
(- \lap)^j w_{\epsilon} > 0 , \;\; j=1, 2, \cdots, k-1.$$
This would contradict with (\ref{4.2}) for sufficiently small $\epsilon$. Hence (\ref{4.2}) must be false and this completes the proof of the theorem.

\section{The Equivalence between Integral and PDE Systems}

In this section, we establish the equivalence between the integral system
\begin{equation}
u_i(x) = \int_{R^n} \frac{1}{|x-y|^{n-\alpha}} f_i(u(y)) d y , \;\; x \in R^n, \; i=1, \cdots, m
\label{ints}
\end{equation}
and the PDE system
\begin{equation}
(- \lap )^{\alpha/2} u_i = f_i(u) , \;\; x \in R^n , \; i = 1, 2, \cdots, m.
\label{PDEs}
\end{equation}
where $u = (u_1, \cdots, u_m)$.

We say that systems (\ref{ints}) and (\ref{PDEs}) are equivalent, if whenever
$u$ is a positive solution of (\ref{ints}), then $u$ is also a solution of
\begin{equation}
(- \lap )^{\alpha/2} u_i =  c f_i(u) , \;\; x \in R^n, \; i=1, \cdots, m
\label{zc}
\end{equation}
with some constant $c$; and vice versa.
\medskip

We first prove

\begin{thm}
Let $\alpha = 2k$ be an even number less than $n$. Then
every positive solution of PDE (\ref{PDEs}) satisfies integral equation (\ref{ints}) with $f_i(u)$ replaced by $c f_i(u)$. \label{thm6.1}
\end{thm}

\begin{rem} i) The converse of Theorem \ref{thm6.1} is much easier to obtain by elementary argument, and we here skip its proof.

ii) The proof of Theorem \ref{mthmequiv2} is entirely similar to that of Theorem \ref{thm6.1}, we also skip it.
\end{rem}

Notice that here we do not assume any asymptotic behavior of the solution near infinity for PDE system (\ref{PDEs}), hence to prove Theorem \ref{thm6.1}, we first need to show that the integral on the right hand side of (\ref{ints}) is finite. We will multiply both sides of (\ref{PDEs}) by proper functions and integrate by parts on $B_r(0)$, then let $r \ra \infty$ and take limits. To this end, we need some estimates of the solutions as stated in the following lemma. One will see that Theorem \ref{thm4.1} (the super polyharmonic property) is a key ingredient in the proof of the following lemma.
\medskip

\begin{lem}
Let $u=(u_1, \cdots, u_m)$ be a positive solution of (\ref{PDEs}). Write
$$F(u) = \sum_{i=1}^m f_i(u); \; \mbox{ and } \; v_{ij} = (- \lap )^j u_i, \;\; j =0, 1, \cdots,
k-1.$$ Then
\be \int_{R^n} \frac{1}{|x|^{n-2k}} F(u(x)) dx \leq C(n) w(0)
< \infty,  \label{6.2} \ee
where $w=u_1 + \cdots + u_m$. We also have
\be \int_{R^n}
\frac{v_{ij}}{|x|^{n-2j}} dx < \infty  \;\; \mbox{ for } j = 1,
\cdots, k-1; \, i=1, \cdots, m .
\label{6.6}
\ee
\label{lem6.2a}
\end{lem}

As an immediate consequence of Lemma \ref{lem6.2a}, we derive
\begin{lem}
There exists a sequence $r_h
\ra \infty$, such that
\be
\int_{\partial B_{r_h}(0)} \sum_{j=0}^{k-1} \frac{v_{ij}}{r_h^{n-2j-1}}  \, d \sigma \ra 0
\; \mbox{ for } i=1, \cdots, m.
\label{A1}
\ee
\label{lem6.2}
\end{lem}

{\bf Proof of Lemma \ref{lem6.2a}}
\medskip

Let $\delta (x)$ be the Dirac Delta function. Let $\phi$ be the
solution of the following boundary value problem
\be
\left\{\begin{array}{ll} (-\lap)^k \phi = \delta (x) & x \in
B_r(0)
\\
\phi = \lap \phi = \cdots = \lap^{k-1} \phi = 0 & \mbox{ on }
\partial B_{r}(0).
\end{array}
\right.
\label{6.7}
\ee

By the maximum principle, one can easily verify that
\be
\frac{\partial}{\partial \nu}[(-\lap)^j \phi] \leq 0, \;\; j = 0,
1, \cdots, k-1, \;\; \mbox{ on }
\partial B_r(0). \label{6.9}
\ee

Sum over the $m$ equations in (\ref{PDEs}), we obtain
\begin{equation}
(- \lap )^k w = F(u(x)) , \;\; x \in R^n .
\label{4.4}
\end{equation}

Multiply both side of the equation (\ref{4.4}) by $\phi$ and
integrate on $B_r(0)$. After integrating by parts several times
and applying Theorem \ref{thm4.1} and (\ref{6.9}), we arrive at
\begin{eqnarray}
\int_{B_r(0)} F(u(x)) \phi (x) dx &=&w(0)+ \sum_{j=0}^{k-1} \int_{\partial
B_r(0)} (-\lap)^j w \,
\frac{\partial}{\partial \nu} [(-\lap)^{k-1-j}\phi] d\sigma \nonumber \\
 &\leq& w(0).
\label{6.8}
\end{eqnarray}

Now letting $r \ra \infty$, one can see that (\ref{6.2}) is just a
direct consequence of the following fact
\be
\phi (x) \ra \frac{c}{|x|^{n-2k}}
\label{6.10}
\ee
with some constant $c$.

To verify (\ref{6.10}), we notice that (\ref{6.7}) is equivalent
to the following system of equations \be \left\{\begin{array}{ll}
- \lap \phi = \psi_1, & \phi
\mid_{\partial B_r} = 0 \\
- \lap \psi_1 = \psi_2, & \psi_1 \mid_{\partial B_r} = 0 \\
\cdots \\
- \lap \psi_{k-1} = \delta (x), & \psi_{k-1} \mid_{\partial B_r} =
0 .
\end{array}
\right. \label{6.11}
\ee
Applying maximum principle consecutively
to $\psi_j$ for $j=1, \cdots, k-1$, one derives that:
$$\psi_{j}(x) \nearrow  c \frac{1}{|x|^{n-2k+2j}} \mbox{ for } j = 0,
1, \cdots, k-1, \mbox{ as } r \ra \infty.$$ This implies
(\ref{6.10}).
\medskip

Now, to derive (\ref{6.6}),  we simply apply the above arguments to
the following equation instead of (\ref{4.4}):
\[
(-\lap)^j u_i = v_{ij},
\]
for each $j=1, \cdots, k-1$ and $i=1, \cdots, m$. This completes the proof of Lemma \ref{lem6.2a}.
\medskip

{\bf The Proof of Lemma \ref{lem6.2}}.

To verify (\ref{A1}) for each fixed $i$, we sum over (\ref{6.2}) and (\ref{6.6}) to obtain
\begin{equation}
\int_{R^n} \left\{ \frac{1}{|x|^{n-2k}} F(u(x)) + \sum_{j=1}^{k-1}
\frac{v_{ij}}{|x|^{n-2j}}\right\} dx < \infty .
\label{z3}
\end{equation}
This implies that there exists a sequence $r_h \ra \infty$, such that
$$
\int_{\partial B_{r_h}(0)} \left\{ \frac{F(u(x))}{r_h^{n-2k-1}} + \sum_{j=1}^{k-1} \frac{v_{ij}}{r_h^{n-2j-1}} \right\} d \sigma \ra 0 .
$$

Since each term in the above summation is nonnegative, we have
\begin{equation}
\frac{1}{r_h^{n-2k-1}} \int_{\partial B_{r_h}} F(u) d \sigma \ra 0
\label{4.5}
\end{equation}
and
\begin{equation}
\int_{\partial B_{r_h}(0)} \sum_{j=1}^{k-1} \frac{v_{ij}}{r_h^{n-2j-1}} d \sigma \ra 0 .
\label{z4}
\end{equation}

For any $\epsilon > 0$, it is obvious that
\begin{equation}
\frac{1}{|\partial B_{r_h}|} \int_{\partial B_{r_h}} w \, d \sigma \leq \epsilon + \frac{1}{|\partial B_{r_h}|}
\int_{\partial B_{r_h}} w \, \chi_{\epsilon} \, d \sigma ,
\label{4.6}
\end{equation}
where $\chi_{\epsilon}$ is the characteristic function on the set
$\{ x \in \partial B_{r_h} \mid w(x) \geq \epsilon \}.$

On the other hand, by condition (\ref{f2}) and Jensen's inequality, we have
\begin{eqnarray*}
 \frac{1}{|\partial B_{r_h}|} \int_{\partial B_{r_h}} F(u) \chi_{\epsilon} \, d \sigma
&\geq& C_{\epsilon} \frac{1}{|\partial B_{r_h}|} \int_{\partial B_{r_h}} ( w \chi_{\epsilon})^p \, d \sigma \\
&\geq& C_{\epsilon} \left( \frac{1}{|\partial B_{r_h}|} \int_{\partial B_{r_h}} w \, \chi_{\epsilon} \, d \sigma
\right)^p .
\end{eqnarray*}
This, together with (\ref{4.5}) and (\ref{4.6}), imply
$$ \frac{1}{|\partial B_{r_h}|} \int_{\partial B_{r_h}} w \, d \sigma \ra 0, \; \mbox{ as } h \ra \infty.$$
Since each $u_i$ is nonnegative, we have
$$ \frac{1}{r_h^{n-1}} \int_{\partial B_{r_h}} u_i d \sigma \ra 0 .$$
Combining this with (\ref{z4}), we arrive at (\ref{A1}).
This completes the proof of Lemma \ref{lem6.2}.
\bigskip

{\bf Proof of Theorem \ref{thm6.1}.}

For each $r > 0$, let $\phi_r(x)$ be the solution of (\ref{6.7}).
Then as in the previous lemma, one verifies that \be \phi_r(x) =
\frac{1}{r^{n-2k}} \phi_1 (\frac{x}{r}) \label{a11}
\ee
and
\be |\phi_1(x)| \leq \frac{C}{|x|^{n-2k}} .
\label{a12}
\ee
It
follows that,
\be \phi_r(x) \leq \frac{C}{|x|^{n-2k}} .
\label{a14}
\ee
Also one can verify that, on $\partial B_r(0)$,
\be |\frac{\partial}{\partial \nu} [(- \lap)^j \phi_r] | \leq
\frac{C}{r^{n-2k +1+2j}} .
\label{a13}
\ee

Multiply both side of the equation
$$ (- \lap)^k u_i = f_i (u) $$ by $\phi_r(x)$ and
integrate on $B_r(0)$. After integrating by parts several times
and applying Theorem \ref{thm4.1} and (\ref{6.9}), we arrive at
\begin{equation}
\int_{B_r(0)} f_i(u(x)) \phi_r (x) dx = u_i(0)+ \sum_{j=0}^{k-1} \int_{\partial
B_r(0)} v_{ij} \,
\frac{\partial}{\partial \nu} [(-\lap)^{k-1-j}\phi_r] d\sigma
\label{6.8a}
\end{equation}

By virtue of (\ref{a13}) and (\ref{A1}),  there exist a sequence $r_h
\ra \infty$, such that the boundary integral on $\partial
B_{r_h}(0)$ in (\ref{6.8a}) approaches 0 as $r_h \ra \infty.$

Applying (\ref{a14}), (\ref{6.2}), and the Lesbegue Dominant Convergence
Theorem to the left hand side of (\ref{6.8a}), and taking limit
along the sequence $\{r_h\}$, we conclude that
$$c \int_{R^n} \frac{1}{|y|^{n-2k}} f_i(u(y)) dy = u_i(0), \;\; i=1, \cdots, m
.$$ By a translation, that is, for each fixed $x$, integrating on $B_r(x)$ instead of on $B_r(0)$, we can derive that $u_i(x)$ is a solution of
\begin{equation}
 u_i(x) = c \int_{R^n} \frac{1}{|x-y|^{n-2k}} f_i(u(y)) d y , \;\; i=1, \cdots, m.
 \label{intsc}
 \end{equation}
This completes the proof of the theorem.
\bigskip

So far, we have proved that if $\alpha =2k$ is an even number, then
every solution of the PDE system (\ref{PDEs}) is a solution of our integral
system (\ref{intsc}).

Now for other real values of $\alpha$, we define the
positive solution of (\ref{PDEs}) in the distribution sense, i.e. $u
\in H^{\frac{\alpha}{2}}(R^n)$ and satisfies
\be
\int_{R^n} (- \lap)^{\frac{\alpha}{4}} u_i \, (- \lap)^{\frac{\alpha}{4}} \phi \, dx =
\int_{R^n} f_i(u(x)) \, \phi (x) \, dx , \;\; i=1, \cdots, m;
\label{6.20}
\ee
for any
$\phi \in C_0^{\infty}$ and $\phi (x) \geq 0.$ Here, as usual,
$$\int_{R^n} (- \lap)^{\frac{\alpha}{4}} u_i \, (-
\lap)^{\frac{\alpha}{4}} \phi \, dx$$
is defined by Fourier transform
$$\int_{R^n} |\xi|^{\alpha} \hat{u_i} (\xi) \overline{\hat{\phi} (\xi)} \,
d\xi , $$ where $\hat{u_i}$ and $\hat{\phi}$ are the Fourier
transform of $u_i$ and $\phi$ respectively.

By taking limits, one can see that (\ref{6.20}) is also true for
any $\phi \in H^{\frac{\alpha}{2}}.$

\begin{thm}
PDE system (\ref{PDEs}) as defined above is
equivalent to integral system (\ref{ints}).
\label{thm6.2}
\end{thm}

{\bf Proof.} (i) For any $\phi \in C_0^{\infty}(R^n)$, let
$$\psi (x) = c \int_{R^n} \frac{\phi (y)}{|x - y|^{n-\alpha}} \, dy
.$$
Choose an appropriate constant $c$, so that $(- \lap)^{\alpha/2} \psi = \phi ,$ consequently $\psi
\in H^{\alpha}(R^n) \subset H^{\frac{\alpha}{2}}(R^n)$, and hence
(\ref{6.20}) holds for $\psi$:
$$\int_{R^n} (- \lap)^{\frac{\alpha}{4}} u_i \, (-
\lap)^{\frac{\alpha}{4}} \psi \, dx = \int_{R^n} f_i(u(x)) \, \psi
(x) \, dx, \;\; i=1, \cdots, m. $$
Integration by parts of the left hand side and
exchange the order of integration of the right hand side yields
$$\int_{R^n} u_i(x) \phi (x) \, dx = \int_{R^n} \left\{ c \int_{R^n}
\frac{f_i(u(y)) }{|x - y|^{n-\alpha}} \, dy \right\} \, \phi (x) \,
dx .$$ Since $\phi$ is any nonnegative $C_0^{\infty}$ function, we
conclude that $u_i$ satisfies the integral equation
$$  u_i (x) = c \int_{R^n} \frac{1}{|x-y|^{n-\alpha}} f_i(u(y)) d y .$$

(ii) Now assume that $u_i \in H^{\frac{\alpha}{2}}(R^n)$ is a
solution of the integral equation (\ref{ints}). Make a Fourier
transform on both sides (cf. [LL], Corollary 5.10),  we have, for some constant $c$,
$$
\hat{u_i}(\xi) = c |\xi|^{-\alpha} \widehat{f_i(u)} (\xi) . $$
 It
follows that $$ \int_{R^n} (- \lap)^{\frac{\alpha}{4}} u_i \, (-
\lap)^{\frac{\alpha}{4}} \phi \, dx = c \int_{R^n}
\widehat{f_i(u)}(\xi) \overline{\hat{\phi}(\xi)} = c \int_{R^n} f_i(u(x))
\phi (x) dx .$$
That is, $u_i$ is a solution of
$$ (- \lap)^{\alpha/2} u_i = c f_i(u(x)) , \;\; x \in R^n $$
in the sense of distributions.

This completes the proof of the theorem.

{\em Authors' Addresses and E-mails:}
\medskip

Wenxiong Chen

Department of Mathematics

Yeshiva University

500 W 185th Street

New York,  NY 10033

wchen@yu.edu
\medskip

Congming Li

Department of Mathematics

University of Colorado at Boulder

Boulder, CO 80309

cli@colorado.edu

\end{document}